\documentclass[12pt]{amsart}

\usepackage{amsmath,amssymb,epsfig,colordvi,graphicx,rotating} 
\newcommand\bfe{{\bf e}} 
\newcommand\bfm{{\bf m}} 
\newcommand\bff{{\bf f}} 
 
\newcommand\x{{\bf x}}

\newcommand{\cc}{{\mathbb C}} 
\newcommand{\pp}{{\mathbb P}} 
 
\newcommand{\nn}{{\mathbb N}}

\copyrightinfo{}{}
\title[Littlewood-Richardson Homotopies]{Solving Schubert Problems with\\
  Littlewood-Richardson Homotopies}

\author{Frank Sottile}
\address[Sottile]{Department of Mathematics\\
         Texas A\&M University\\
         College Station\\
         TX \ 77843}
\email{sottile@math.tamu.edu}
\urladdr{http://www.math.tamu.edu/\~{}sottile/}
\thanks{Work of Sottile supported by the NSF
under Grants DMS-0538734 and DMS-0915211.}

\author{Ravi Vakil}
\address[Vakil]{Department of Mathematics\\
      Stanford University\\ 
       Stanford, CA 94305}
\email{vakil@math.stanford.edu}
\urladdr{ http://math.stanford.edu/{\~{}}vakil}
\thanks{Work of Vakil supported by
the NSF under Grants DMS-0801196.}

\author{Jan Verschelde}
\address[Verschelde]{Department of Mathematics, Statistics,
and Computer Science\\
University of Illinois at Chicago\\
851 South Morgan (M/C 249)\\
Chicago, IL 60607-7045}
\email{jan@math.uic.edu}
\urladdr{http://www.math.uic.edu/{\~{}}jan}
\thanks{Work of Verschelde supported by the NSF
under Grant No.\ 0713018.}

\subjclass[2000]{Primary 65H10.  Secondary 14N15, 14Q99, 68W30} 
 
\keywords{ continuation, geometric Littlewood-Richard\-son rule, Grassmannian, homotopies,
  numerical Schubert calculus, path following, polynomial system, Schubert problems.}
 
\begin{document} 

\begin{abstract} 
We present a new numerical homotopy continuation algorithm for finding  
all solutions to Schubert problems on Grassmannians. 
This Littlewood-Richardson homotopy is based on Vakil's geometric proof 
of the Littlewood-Richardson rule. 
Its start solutions are given by linear equations and they are tracked  
through a sequence of homotopies encoded by certain checker configurations  
to find the solutions to a given Schubert problem. 
For generic Schubert problems the number of paths tracked is optimal. 
The Littlewood-Richardson homotopy algorithm is implemented 
using the path trackers of the software package PHCpack. \vspace{-20pt}
\end{abstract} 
%%%%%%%%%%%%%%%%%%%%%%%%%%%%%%%%%%%%%%%%%%%%%%%%%%%%%%%%%%%%%%%%%%%%%%%%%%% 

\maketitle 
  
%%%%%%%%%%%%%%%%%%%%%%%%%%%%%%%%%%%%%%%%%%%%%%%%%%%%%%%%%%%%%%%%%%%%%%%%%%% 
 
\section{Introduction} 
The Schubert calculus is concerned with geometric problems of the form: 
Determine the $k$-dimensional linear subspaces of $\cc^n$ that meet  
a collection of fixed linear subspaces in specified dimensions. 
For example, what are the three-dimensional linear subspaces of $\cc^7$  
that meet each of 12 general four-dimensional linear subspaces in at  
least a line?  (There are $462$~\cite{Sch1886}.) 
The traditional goal is to count the number of solutions and the method  
of choice for this enumeration is the Littlewood-Richardson rule,  
which comes from combinatorics and representation theory~\cite{Ful97}. 
Recently, Vakil gave a geometric proof of this rule~\cite{Vak06} through  
explicit specializations organized by a combinatorial checkers game. 
 
Interest has grown in  computing the solutions to actual Schubert problems. 
One motivation has been the experimental study of reality in the Schubert 
calculus~\cite{HGMRTJS09,Sot03,Sot10,Sot00}. 
A  proof of Pieri's rule (a special case of the Littlewood-Richardson 
rule) using geometric specializations~\cite{Sot97a} led to the Pieri  
homotopy for solving special Schubert problems~\cite{HSS98}. 
This was implemented and refined~\cite{HV00,LWW02,Ver00,VW04,VW04b}, 
and has been used to address a problem in pure mathematics~\cite{LS09}. 
Another motivation is the output pole placement problem
in linear systems control~\cite{BB81,Byr89,EG02,RRW96,VW04}.
 
We present the Littlewood-Richardson homotopy,  
which is a numerical homotopy algorithm 
for finding all solutions to any Schubert problem. 
It is based on the geometric Littlewood-Richardson rule~\cite{Vak06}  
and it is optimal in that generically there are no extraneous paths  
to be tracked. 

We describe Schubert problems and their equations in \S\ref{S:schubert},
and give a detailed example of the geometric Littlewood-Richardson rule in
\S\ref{fourlines}.
We then explain the local structure of the Littlewood-Richardson homotopy in \S\ref{six}.
The next three sections give more details on the local coordinates, the moving flag, and
the checker configurations.
In \S\ref{S:solving}  we discuss the global structure of the Littlewood-Richardson
homotopy and conclude in \S\ref{S:computation} with a brief description of our
PHCpack~\cite{Ver99} implementation and timings.

%%%%%%%%%%%%%%%%%%%%%%%%%%%%%%%%%%%%%%%%%%%%%%%%%%%%%%%%%%%%%%%%%%%%%%%%%%% 
% 
\section{Schubert Problems}\label{S:schubert} 
 
A Schubert problem asks for the $k$-dimensional subspaces of $\cc^n$  
that satisfy certain Schubert conditions imposed by general flags. 
We explain this in concrete terms. 
 
A point in $\cc^n$ is represented by a $n\times 1$ column vector and  
a linear subspace as the column span of a matrix. 
A flag $F$ is represented by an ordered 
basis $\bff_1, \dotsc, \bff_n $ of $\cc^n$ that 
forms the columns of a matrix $F$. 
If we write $F_i$ for the span of the first $i$ columns of $F$, 
then a Schubert condition  
imposed by $F$ is the condition on the $k$-plane $X$ that 
 \begin{equation}\label{Eq:Schubert_condition} 
  \dim ( X\cap F_{\omega_i} )\geq i 
   \quad\mbox{for}\quad 1\leq i\leq k\,, 
 \end{equation} 
where $\omega\in\nn^k$ is a bracket; 
$1 \leq \omega_1 < \omega_2 < \cdots < \omega_k \leq n$. 
 
If we set $|\omega|=\sum_i n-k+i-w_i$, then  
a Schubert problem is a list 
$\omega^1,\omega^2,\dotsc,\omega^s$ of brackets such that  
 \begin{equation}\label{Eq:codim_condition} 
   |\omega^1|+|\omega^2|+\dotsb+|\omega^s|\ =\ k(n-k)\,. 
 \end{equation} 
For example, the Schubert problem of three-planes meeting 12 four-planes  
in $\cc^7$ is given by 12 equal codimension one brackets and is written 
succinctly as $[4~6~7]^{12}$.  
The  numerical condition~\eqref{Eq:codim_condition}  
ensures that if $F^1,\dotsc,F^s$ 
are general, then there are finitely many $k$-planes that satisfy  
condition $\omega^i$ 
for flag $F^i$, for $i=1,\dotsc,s$. 
 
The set of $k$-planes $X$ satisfying~\eqref{Eq:Schubert_condition} 
is the Schubert variety $\Omega_\omega(F)$. 
This is a subvariety of the $k(n{-}k)$-dimensional Grassmannian of $k$-planes 
in $n$-space. 
Thus solving a Schubert problem corresponds to determining the intersection of 
Schubert varieties with respect to various flags. 
 
These geometric conditions are formulated as systems of polynomials 
by parameterizing an appropriate subset of the Grassmannian. 
For example, for $F \in \cc^{6 \times 6}$, the Schubert variety $\Omega_{[2~4~6]}(F)$  
contains 
\begin{equation}\label{Eq:rank_condition} 
   X = \left[ 
     \begin{array}{ccc} 
        1 & 0 & 0 \\ 
        x_{21} & 1 & 0 \\  
        x_{31} & x_{32} & 1 \\ 
        x_{41} & x_{42} & x_{43} \\ 
         0 & x_{52} & x_{53} \\ 
         0 & 0 & x_{63} \\ 
     \end{array} 
   \right] 
   \quad 
   \begin{array}{l} 
      \dim(X \cap F_2 ) = 1 \\ 
      \dim(X \cap  F_4 ) = 2 \\ 
      \dim(X \cap F_6 ) = 3  
   \end{array} 
\end{equation} 
Expressed via conditions on the minors of $[X|F_i]$ this is a system of 13 
polynomials in 9 variables. 
 
The most elementary Schubert problem involves only two brackets, $\omega$ and $\tau$ with 
$|\omega|+|\tau|=k(n{-}k)$.  
If $F$ and $M$ are general flags and  
$\omega^\vee=[n{+}1{-}\omega_k\;\dotsc~n{+}1{-}\;\omega_1]$, then  
 \begin{equation}\label{Eq:double_intersection} 
  \Omega_\omega(F)\cap\Omega_\tau(M)= 
   \left\{\begin{array}{cl}  
   \langle {\bf x}_1,\dotsc,{\bf x}_k \rangle&\mbox{if }\tau=\omega^\vee,\\ 
   \emptyset&\mbox{otherwise,}\end{array}\right.  
 \end{equation} 
where  
${\bf x}_i=F_{\omega_i}\cap M_{n+1-\omega_i}$, which is  
one-dimensional and thus solved by linear algebra.
%Thus this is solved by linear algebra. 
Such elementary Schubert problems are the 
start systems for the Littlewood-Richardson homotopy.  
 
When $|\omega|+|\tau|<k(n{-}k)$, an intersection  
 \begin{equation}\label{Eq:general_double} 
    \Omega_\omega(F)\cap\Omega_\tau(M) 
 \end{equation} 
of Schubert varieties for general flags $F$ and $M$ has positive dimension. 
This intersection is homologous to  
a union of Schubert varieties $\Omega_\sigma(F)$ for  
$|\sigma|=|\omega|+|\tau|$, each occurring  
with multiplicity the Littlewood-Richardson number 
$c^\sigma_{\omega,\tau}$. 
We write this formally as a sum, 
\begin{equation} 
   \Omega_\omega(F)\cap\Omega_\tau(M) \sim
   \sum_\sigma c^\sigma_{\omega,\tau} \Omega_\sigma(F). 
\end{equation} 
 
In the geometric Littlewood-Richardson rule~\cite{Vak06},  
the flag $M$ moves into special 
position with respect to the flag $F$. 
This changes the intersection~\eqref{Eq:general_double},  
breaking it into components which are transformed into Schubert 
varieties $\Omega_\sigma(F)$. 
Then $c^\sigma_{\omega,\tau}$ is the number of different ways to arrive at  
$\Omega_\sigma(F)$. 
 
\begin{figure}[hbt] 
\[
  \begin{picture}(76,66)(-19,-1) 
    \put(-3,-3){\includegraphics{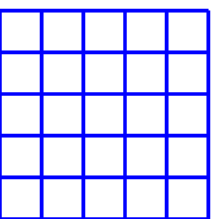}} 
    \put( 0,0){\small$1$}\put( 0,12){\small$1$}\put( 0,24){\small$1$}
            \put( 0,36){\Red{\small$1$}}\put( 0,48){\small$0$} 
    \put(12,0){\small$2$}\put(12,12){\Red{\small$2$}}\put(12,24){\small$1$}
           \put(12,36){\small$1$}\put(12,48){\small$0$} 
    \put(24,0){\Red{\small$3$}}\put(24,12){\small$2$}\put(24,24){\small$1$}
           \put(24,36){\small$1$}\put(24,48){\small$0$} 
    \put(36,0){\small$4$}\put(36,12){\small$3$}\put(36,24){\Red{\small$2$}}
           \put(36,36){\small$1$}\put(36,48){\small$0$} 
    \put(48,0){\small$5$}\put(48,12){\small$4$}\put(48,24){\small$3$}
           \put(48,36){\small$2$}\put(48,48){\Red{\small$1$}} 
  
    \put(-19,1){\small$M_5$}\put(-19,13){\small$M_4$}\put(-19,25){\small$M_3$}
           \put(-19,37){\small$M_2$}\put(-19,49){\small$M_1$} 
    \put(-2,61){\small$F_1$}\put(10,61){\small$F_2$}\put(22,61){\small$F_3$}
           \put(34,61){\small$F_4$}\put(46,61){\small$F_5$}  
  \end{picture}
   \quad\raisebox{27pt}{$\Longleftrightarrow$}\quad 
  \begin{picture}(70,66)(0,2) 
    \put(0,0){\includegraphics{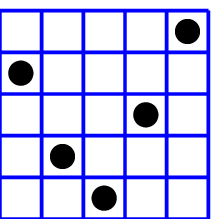}} 
  \end{picture}\vspace{-3pt}
\]
\caption{Dimension array $\dim M_a \cap F_b$ 
         and corresponding checker configuration.}\vspace{-3pt}
\label{figdimarray} 
\end{figure}

The relative position of the flags $F$ and $M$ is represented  
via a configuration of $n$ black 
checkers in a $n\times n$ board with no two in the same row or column. 
The dimension of $M_a\cap F_b$ is the number of checkers weakly northwest  
of the square $(a,b)$.   
This is illustrated in Figure~\ref{figdimarray}. 
Each cell corresponds to a vector space, and the vector  
space of each cell is contains the vector spaces of the cells weakly northwest 
of it. 
 
\begin{figure}[hbt] 
\[
   \includegraphics{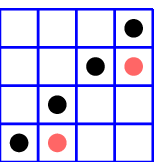}\qquad 
   \includegraphics{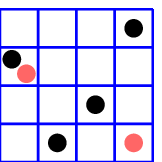}\qquad 
   \includegraphics{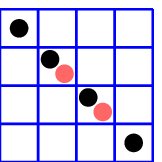} \vspace{-3pt}
\]
\caption{Three checkerboards with $n = 4$ and $k = 2$.} \vspace{-3pt}
\label{figthreeboards} 
\end{figure}

All components of the specializations of~\eqref{Eq:general_double}  
are represented by placements  
of $k$ red\footnote{Red checkers look grey when printed in
black and white.}
checkers on a board 
with $n$ checkers representing the relative positions of 
the flags. 
The red checkers represent the position of a typical $k$-plane in the 
component as follows:  If the $k$-plane meets the vector space 
corresponding to a cell in dimension $\ell$, then there are $\ell$ red 
checkers weakly northwest of it. 
See Figure~\ref{figthreeboards}
 for examples. 
We discuss the placement and movement of the checkers  
in \S \ref{fourlines} and \S \ref{checkers}. 
 
Applying the geometric Littlewood-Richardson rule to two Schubert  
varieties in a Schubert problem of $s$ brackets reduces it to Schubert  
problems involving $s{-}1$ brackets.  
The Littlewood-Richardson homotopy begins with the solutions to  
those smaller problems and 
reverses the specializations to solve the original Schubert problem. 
 
Let $k=3$ and $n=6$, and consider this for the Schubert problem  
$[2~4~6]^3 = [2~4~6] [2~4~6] [2~4~6]$. 
Given three general flags $F,M,N$, we want to resolve the triple intersection 
 \begin{equation}\label{Eq:triple} 
   \Omega_{[2~4~6]}(F) \cap \Omega_{[2~4~6]}(M) \cap \Omega_{[2~4~6]}(N). 
 \end{equation} 
We first apply the geometric Littlewood-Richardson rule to the first  
intersection to obtain 
 \begin{equation}\label{Eq:GLRR} 
   \bigl(\Omega_{[2~3~4]}(F)+2\Omega_{[1~3~5]}(F)+\Omega_{[1~2~6]}(F)\bigr) \cap \Omega_{[2~4~6]}(N), 
 \end{equation} 
and then apply~\eqref{Eq:double_intersection} to obtain 
%% 
% \begin{equation}\label{Eq:solution} 
$    2\langle {\bf x}_1,{\bf x}_2,{\bf x}_3 \rangle, $
% \end{equation} 
% 
where ${\bf x}_1=F_1\cap N_6$,  
${\bf x}_2=F_3\cap N_4$, and ${\bf x}_3=F_5\cap N_2$. 
 
The Littlewood-Richardson homotopy starts with the single 
3-plane $\langle {\bf x}_1,{\bf x}_2,{\bf x}_3 \rangle$ (counted twice) which is the
unique solution  
to~\eqref{Eq:GLRR}.  
It then numerically continues this solution backwards along the geometric specializations 
transforming~\eqref{Eq:triple} into~\eqref{Eq:GLRR} to arrive at solutions 
to~\eqref{Eq:triple}. 
As the multiplicity 2 of $\Omega_{[1~3~5]}(F)$ in~\eqref{Eq:GLRR} is the number of paths in 
the specialization that end in $\Omega_{[1~3~5]}(F)$, the single solution 
$\langle {\bf x}_1,{\bf x}_2,{\bf x}_3 \rangle$ that we began with  
yields two solutions to~\eqref{Eq:triple}.

%%%%%%%%%%%%%%%%%%%%%%%%%%%%%%%%%%%%%%%%%%%%%%%%%%%%%%%%%%%%%%%%%%%%%%%%%%%%%%%%%%%%%%%% 
% 
\section{The problem of four lines} 
\label{fourlines} 
 
We illustrate the Littlewood-Richardson homotopy via the 
classical problem of which lines in projective three-space ($\pp^3$) meet  
four given lines. 
This corresponds to two-planes in $\cc^4$ meeting four fixed two-planes  
nontrivially, or $[2~4]^4$.

%%%%%%%%%%%%%%%%%%%%%%%%%%%%%%%%%%%%%%%%%%%%%%%%%%%%%%%%%%%%%%%%%%%%% 
\begin{figure*}[hbt] 
\[ 
  \begin{picture}(470,188)(0,13) 
   \put(  0,100){\includegraphics[height=3.5cm]{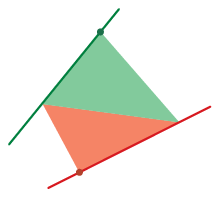}} 
   \put(125,100){\includegraphics[height=3.5cm]{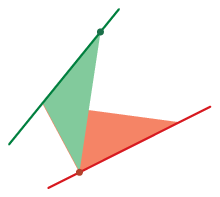}} 
   \put(250,100){\includegraphics[height=3.5cm]{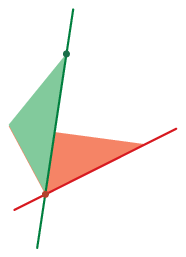}} 
   \put(375,100){\includegraphics[height=3.5cm]{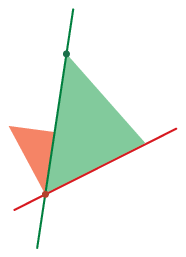}} 
   \put( 56,  5){\includegraphics[height=3.5cm]{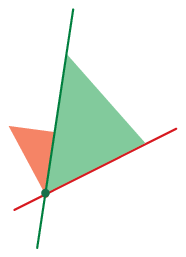}} 
   \put(180,  5){\includegraphics[height=3.5cm]{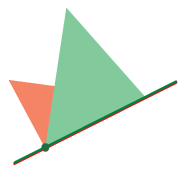}} 
   \put(304,  5){\includegraphics[height=3.5cm]{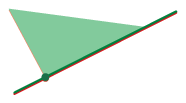}} 
   \put( 10,190){stage 0} 
   \put(135,190){stage 1}    \put( 65,95){stage 4} 
   \put(260,190){stage 2}    \put(190,95){stage 5} 
   \put(385,190){stage 3}    \put(315,95){stage 6} 
   \put(110,145){\Large$\rightarrow$} \put( 45,50){\Large$\rightarrow$} 
   \put(235,145){\Large$\rightarrow$} \put(170,50){\Large$\rightarrow$} 
   \put(360,145){\Large$\rightarrow$} \put(295,50){\Large$\rightarrow$} 
  \end{picture} 
\] 
\[ 
  \begin{picture}(432,83)(0,8) 
   \put(0,0){ 
    \begin{picture}(36,91)(-2,-43) 
     \put(1,41){stage 0} 
     \put(0,27){\small$1$} \put(9,27){\small$1$} \put(18,27){\small$1$} \put(27,27){\small$1$} 
     \put(0,18){\small$1$} \put(9,18){\small$1$} \put(18,18){\small$1$} \put(27,18){\small$0$} 
     \put(0, 9){\small$1$} \put(9, 9){\small$1$} \put(18, 9){\small$0$} \put(27, 9){\small$0$} 
     \put(0, 0){\small$1$} \put(9, 0){\small$0$} \put(18, 0){\small$0$} \put(27, 0){\small$0$} 
     \put(-2,-43){\includegraphics{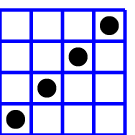}} 
%     \put(2,-54){4 3 2 1} 
    \end{picture}} 
   \put(48,48){\Large$\rightarrow$} 
   \put(66,0){ 
    \begin{picture}(36,91)(-2,-43) 
     \put(1,41){stage 1} 
     \put(0,27){\small$1$} \put(9,27){\small$1$} \put(18,27){\small$1$} \put(27,27){\small$0$} 
     \put(0,18){\small$1$} \put(9,18){\small$1$} \put(18,18){\small$0$} \put(27,18){\small$1$} 
     \put(0, 9){\small$1$} \put(9, 9){\small$1$} \put(18, 9){\small$0$} \put(27, 9){\small$0$} 
     \put(0, 0){\small$1$} \put(9, 0){\small$0$} \put(18, 0){\small$0$} \put(27, 0){\small$0$} 
     \put(-2,-43){\includegraphics{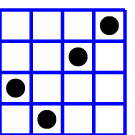}} 
%     \put(2,-54){4 3 1 2} 
    \end{picture}} 
   \put(114,48){\Large$\rightarrow$} 
   \put(132,0){ 
    \begin{picture}(36,91)(-2,-43) 
     \put(1,41){stage 2} 
     \put(0,27){\small$1$} \put(9,27){\small$1$} \put(18,27){\small$0$} \put(27,27){\small$0$} 
     \put(0,18){\small$1$} \put(9,18){\small$0$} \put(18,18){\small$1$} \put(27,18){\small$1$} 
     \put(0, 9){\small$1$} \put(9, 9){\small$0$} \put(18, 9){\small$1$} \put(27, 9){\small$0$} 
     \put(0, 0){\small$1$} \put(9, 0){\small$0$} \put(18, 0){\small$0$} \put(27, 0){\small$0$} 
     \put(-2,-43){\includegraphics{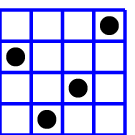}} 
%     \put(2,-54){4 1 3 2} 
    \end{picture}} 
   \put(180,48){\Large$\rightarrow$} 
   \put(198,0){ 
    \begin{picture}(36,91)(-2,-43) 
     \put(1,41){stage 3} 
     \put(0,27){\small$1$} \put(9,27){\small$1$} \put(18,27){\small$0$} \put(27,27){\small$0$} 
     \put(0,18){\small$1$} \put(9,18){\small$0$} \put(18,18){\small$1$} \put(27,18){\small$0$} 
     \put(0, 9){\small$1$} \put(9, 9){\small$0$} \put(18, 9){\small$0$} \put(27, 9){\small$1$} 
     \put(0, 0){\small$1$} \put(9, 0){\small$0$} \put(18, 0){\small$0$} \put(27, 0){\small$0$} 
     \put(-2,-43){\includegraphics{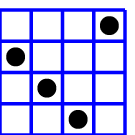}} 
%     \put(2,-54){4 1 2 3} 
    \end{picture}} 
   \put(246,48){\Large$\rightarrow$} 
   \put(264,0){ 
    \begin{picture}(36,91)(-2,-43) 
     \put(1,41){stage 4} 
     \put(0,27){\small$1$} \put(9,27){\small$0$} \put(18,27){\small$0$} \put(27,27){\small$0$} 
     \put(0,18){\small$0$} \put(9,18){\small$1$} \put(18,18){\small$1$} \put(27,18){\small$0$} 
     \put(0, 9){\small$0$} \put(9, 9){\small$1$} \put(18, 9){\small$0$} \put(27, 9){\small$1$} 
     \put(0, 0){\small$0$} \put(9, 0){\small$1$} \put(18, 0){\small$0$} \put(27, 0){\small$0$} 
     \put(-2,-43){\includegraphics{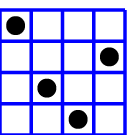}} 
%     \put(2,-54){1 4 2 3} 
    \end{picture}} 
   \put(312,48){\Large$\rightarrow$} 
   \put(330,0){ 
    \begin{picture}(36,91)(-2,-43) 
     \put(1,41){stage 5} 
     \put(0,27){\small$1$} \put(9,27){\small$0$} \put(18,27){\small$0$} \put(27,27){\small$0$} 
     \put(0,18){\small$0$} \put(9,18){\small$1$} \put(18,18){\small$0$} \put(27,18){\small$0$} 
     \put(0, 9){\small$0$} \put(9, 9){\small$0$} \put(18, 9){\small$1$} \put(27, 9){\small$1$} 
     \put(0, 0){\small$0$} \put(9, 0){\small$0$} \put(18, 0){\small$1$} \put(27, 0){\small$0$} 
     \put(-2,-43){\includegraphics{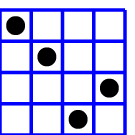}} 
%     \put(2,-54){1 2 4 3} 
    \end{picture}} 
   \put(378,48){\Large$\rightarrow$} 
   \put(396,0){ 
    \begin{picture}(36,91)(-2,-43) 
     \put(1,41){stage 6} 
     \put(0,27){\small$1$} \put(9,27){\small$0$} \put(18,27){\small$0$} \put(27,27){\small$0$} 
     \put(0,18){\small$0$} \put(9,18){\small$1$} \put(18,18){\small$0$} \put(27,18){\small$0$} 
     \put(0, 9){\small$0$} \put(9, 9){\small$0$} \put(18, 9){\small$1$} \put(27, 9){\small$0$} 
     \put(0, 0){\small$0$} \put(9, 0){\small$0$} \put(18, 0){\small$0$} \put(27, 0){\small$1$} 
     \put(-2,-43){\includegraphics{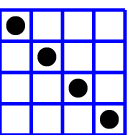}} 
%     \put(2,-54){1 2 3 4} 
    \end{picture}}  
  \end{picture} 
\] 
\caption{Specialization of the moving flag to the fixed flag.}\vspace{-1pt}
\label{figmovingflag} 
\end{figure*} 
%%%%%%%%%%%%%%%%%%%%%%%%%%%%%%%%%%%%%%%%%%%%%%%%%%%%%%%%%%%%%%%%%%%% 
A flag in $\pp^3$ consists of a point lying on a line that is contained  
in a plane (depicted here as a triangle). 
Figure~\ref{figmovingflag} shows the specialization of two flags, one 
fixed and one moving, that underlies the geometric Littlewood-Richardson  
rule for every Schubert problem in $\pp^3$.  
The  top shows the geometry of the specialization. 
Below are  matrices representing the moving flag and checkerboards representing 
the relative positions of the two flags. % and the corresponding permutations. 
We recognize this as the bubble sort of the black checkers.
 
From stage 0 to stage 1 only the moving plane moves; 
the line and point are fixed.   The plane moves until it contains the 
fixed point of the fixed flag.  
From stage $1$ to stage $2$ 
only the moving line moves, until it too contains the fixed point.  
Then the moving plane moves again (to contain the fixed line); then the moving point; then the 
moving line; then the moving plane.  At the end the two flags coincide.
 
We describe the solution to the problem of four lines using  
the Littlewood-Richardson homotopy. 
Let $\ell_1,\ell_2,\ell_3,\ell_4$ be the four lines, where 
$\ell_1$ is the line of the fixed flag, $\ell_2$ the line 
of the moving flag, and $\ell_3$ and $\ell_4$ are two other general lines.   
The family of lines meeting $\ell_1$ and $\ell_2$ is two-dimensional; 
it is parameterized by $\ell_1 \times \ell_2$, as any line   
meeting $\ell_1$ and $\ell_2$ 
is determined by the points where it meets them. 
On this parameterized surface, we seek those points corresponding to 
lines satisfying the further condition $(*)$ of meeting $\ell_3$ and 
$\ell_4$.  Between stages 0 and 1, nothing changes, but in moving to stage  
2, $\ell_2$ moves to intersect $\ell_1$.  There are now two distinct 
two-dimensional families of lines meeting both $\ell_1$ and $\ell_2$: 
(a) those lines lying in the plane $P$ containing $\ell_1$ 
and $\ell_2$ %(or equivalently, those lines meeting $\ell_1$ and 
%$\ell_2$ but not $\ell_1 \cap \ell_2$), 
and (b) those lines in space passing through 
the point $p=\ell_1 \cap \ell_2$.  We now impose the 
additional condition $(*)$ on both of these cases.  In case (a), 
$\ell_i$ meets $P$ in a point $p_i$ ($i = 3,4$), so there is one line  in 
$P$ meeting $\ell_3$ and $\ell_4$, namely $\overline{ p_3 p_4 }$.  In 
case (b), there is one line through $p$ meeting $\ell_3$ and $\ell_4$,
namely $\overline{p,\ell_3}\cap\overline{p,\ell_4}$.
%: standing at the point $p$ and looking outwards, lines $\ell_3$ and $\ell_4$ 
%appear to intersect, and that line from $p$ to the apparent 
%intersection will be a line meeting $\ell_1,\dotsc,\ell_4$.   
After this, the only change is that the plane $P$,  
which equals the moving plane after 
stage 3, and rotates into the fixed plane between stages 5 and 6.    
To solve the original problem, we reverse this process, 
starting with the two solutions in cases (a) and (b), and reversing 
the specialization.   
 
Note that we have reduced one problem involving $4$ brackets to two 
problems involving $3$ brackets. 

\begin{figure*}[hbt] 
\[ 
  \begin{picture}(470,200)(0,8) 
   \put(  0,100){\includegraphics[height=3.5cm]{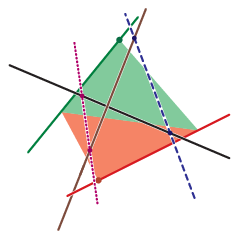}} 
   \put(125,100){\includegraphics[height=3.5cm]{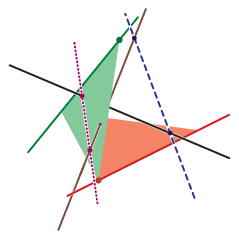}} 
   \put(250,100){\includegraphics[height=3.5cm]{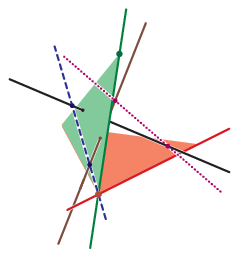}} 
   \put(375,100){\includegraphics[height=3.5cm]{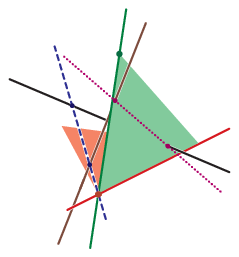}} 
   \put( 56,  0){\includegraphics[height=3.5cm]{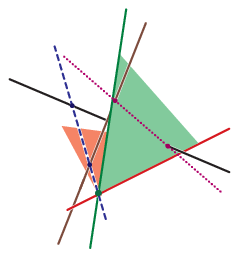}} 
   \put(180,  0){\includegraphics[height=3.5cm]{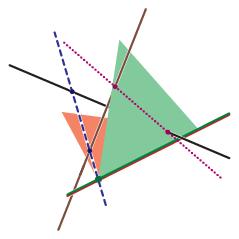}} 
   \put(304,  0){\includegraphics[height=3.5cm]{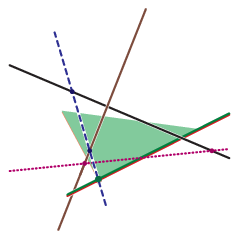}} 
   \put( 10,195){stage 0} 
   \put(135,195){stage 1}    \put( 65,95){stage 4} 
   \put(260,195){stage 2}    \put(190,95){stage 5} 
   \put(385,195){stage 3}    \put(315,95){stage 6} 
   \put(110,145){\Large$\rightarrow$} \put( 45,50){\Large$\rightarrow$} 
   \put(235,145){\Large$\rightarrow$} \put(170,50){\Large$\rightarrow$} 
   \put(360,145){\Large$\rightarrow$} \put(295,50){\Large$\rightarrow$} 
  \end{picture} 
\] 
\[ 
  \begin{picture}(458,195)(0,12) 
% 
%   Make figure left-to-right 
% 
   \put(  8,194){stage 0} 
   \put( 77,194){stage 1} 
   \put(146,194){stage 2} 
   \put(215,194){stage 3} 
   \put(284,194){stage 4} 
   \put(353,194){stage 5} 
   \put(422,194){stage 6} 
  \put(12,125){ 
    \begin{picture}(36,36) 
     \put(0,27){$*$} \put(9,27){\small$0$} 
     \put(0,18){\small$1$}   \put(9,18){\small$0$} 
     \put(0, 9){\small$0$}   \put(9, 9){$*$} 
     \put(0, 0){\small$0$}   \put(9, 0){\small$1$} 
    \end{picture}} 
  \put(0,74){\includegraphics{stage0.eps}} 
  \put(12,35){ 
    \begin{picture}(36,36) 
     \put(0,27){$*$} \put(9,27){\small$0$} 
     \put(0,18){\small$1$}   \put(9,18){\small$0$} 
     \put(0, 9){\small$0$}   \put(9, 9){$*$} 
     \put(0, 0){\small$0$}   \put(9, 0){\small$1$} 
   \end{picture}} 
 
  \put(49,96){\vector(1,0){15}} 
 
  \put(81,125){ 
    \begin{picture}(36,36) 
     \put(0,27){$*$} \put(9,27){\small$0$} 
     \put(0,18){\small$1$}   \put(9,18){\small$0$} 
     \put(0, 9){\small$0$}   \put(9, 9){$*$} 
     \put(0, 0){\small$0$}   \put(9, 0){\small$1$} 
    \end{picture}} 
  \put(69,74){\includegraphics{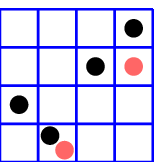}} 
  \put(81,35){ 
    \begin{picture}(36,36) 
     \put(0,27){$*$} \put(9,27){\small$0$} 
     \put(0,18){\small$1$}   \put(9,18){\small$0$} 
     \put(0, 9){\small$0$}   \put(9, 9){$*$} 
     \put(0, 0){\small$0$}   \put(9, 0){\small$1$} 
   \end{picture}} 
 
  \put(118,100){\vector(2, 3){15}} 
  \put(118, 92){\vector(2,-3){15}} 
 
  \put(150,153){ 
    \begin{picture}(36,36) 
     \put(0,27){\small$0$}   \put(9,27){$*$} 
     \put(0,18){\small$1$}   \put(9,18){\small$0$} 
     \put(0, 9){\small$0$}   \put(9, 9){$*$} 
     \put(0, 0){\small$0$}   \put(9, 0){\small$1$} 
    \end{picture}} 
  \put(138,102){\includegraphics{stage2.1.eps}} 
  \put(138, 46){\includegraphics{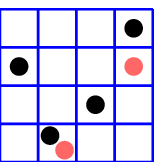}} 
  \put(150,7){ 
    \begin{picture}(36,36) 
     \put(0,27){$*$} \put(9,27){\small$0$} 
     \put(0,18){\small$1$}   \put(9,18){$*$} 
     \put(0, 9){\small$0$}   \put(9, 9){\small$0$} 
     \put(0, 0){\small$0$}   \put(9, 0){\small$1$} 
   \end{picture}} 
 
  \put(187,124){\vector(1,0){15}} 
  \put(187, 68){\vector(1,0){15}} 
 
   \put(219,153){ 
    \begin{picture}(36,36) 
     \put(0,27){\small$0$}   \put(9,27){$*$} 
     \put(0,18){\small$1$}   \put(9,18){\small$0$} 
     \put(0, 9){\small$0$}   \put(9, 9){$*$} 
     \put(0, 0){\small$0$}   \put(9, 0){\small$1$} 
    \end{picture}} 
  \put(207,102){\includegraphics{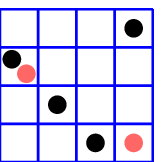}} 
  \put(207, 46){\includegraphics{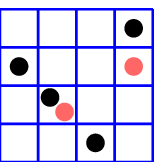}} 
  \put(219,7){ 
    \begin{picture}(36,36) 
     \put(0,27){$*$} \put(9,27){\small$0$} 
     \put(0,18){\small$1$}   \put(9,18){$*$} 
     \put(0, 9){\small$0$}   \put(9, 9){\small$1$} 
     \put(0, 0){\small$0$}   \put(9, 0){\small$0$} 
   \end{picture}} 
 
  \put(256,124){\vector(1,0){15}} 
  \put(256, 68){\vector(1,0){15}} 
 
 \put(288,153){ 
    \begin{picture}(36,36) 
     \put(0,27){\small$1$}   \put(9,27){\small$0$} 
     \put(0,18){\small$0$}   \put(9,18){$*$} 
     \put(0, 9){\small$0$}   \put(9, 9){$*$} 
     \put(0, 0){\small$0$}   \put(9, 0){\small$1$} 
    \end{picture}} 
  \put(276,102){\includegraphics{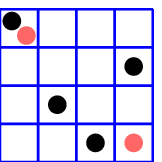}} 
  \put(276, 46){\includegraphics{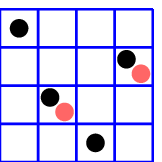}} 
  \put(288,7){ 
    \begin{picture}(36,36) 
     \put(0,27){$*$} \put(9,27){$*$} 
     \put(0,18){\small$1$}   \put(9,18){\small$0$} 
     \put(0, 9){\small$0$}   \put(9, 9){\small$1$} 
     \put(0, 0){\small$0$}   \put(9, 0){\small$0$} 
   \end{picture}} 
 
  \put(325,124){\vector(1,0){15}} 
  \put(325, 68){\vector(1,0){15}} 
 
 \put(357,153){ 
    \begin{picture}(36,36) 
     \put(0,27){\small$1$}   \put(9,27){\small$0$} 
     \put(0,18){\small$0$}   \put(9,18){$*$} 
     \put(0, 9){\small$0$}   \put(9, 9){$*$} 
     \put(0, 0){\small$0$}   \put(9, 0){\small$1$} 
    \end{picture}} 
  \put(345,102){\includegraphics{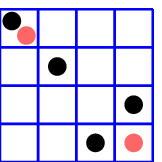}} 
  \put(345, 46){\includegraphics{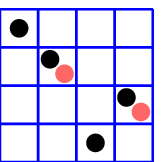}} 
  \put(357,7){ 
    \begin{picture}(36,36) 
     \put(0,27){$*$} \put(9,27){$*$} 
     \put(0,18){\small$1$}   \put(9,18){\small$0$} 
     \put(0, 9){\small$0$}   \put(9, 9){\small$1$} 
     \put(0, 0){\small$0$}   \put(9, 0){\small$0$} 
   \end{picture}} 
 
  \put(394,124){\vector(1,0){15}} 
  \put(394, 68){\vector(1,0){15}} 
 
 \put(426,153){ 
    \begin{picture}(36,36) 
     \put(0,27){\small$1$}   \put(9,27){\small$0$} 
     \put(0,18){\small$0$}   \put(9,18){$*$} 
     \put(0, 9){\small$0$}   \put(9, 9){$*$} 
     \put(0, 0){\small$0$}   \put(9, 0){\small$1$} 
    \end{picture}} 
  \put(414,102){\includegraphics{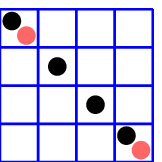}} 
  \put(414, 46){\includegraphics{stage6.2.eps}} 
  \put(426,7){ 
    \begin{picture}(36,36) 
     \put(0,27){$*$} \put(9,27){$*$} 
     \put(0,18){\small$1$}   \put(9,18){\small$0$} 
     \put(0, 9){\small$0$}   \put(9, 9){\small$1$} 
     \put(0, 0){\small$0$}   \put(9, 0){\small$0$} 
   \end{picture}} 
 
  \end{picture} 
\] 
\caption{Resolving the problem of four lines.} 
\label{figcheckergame} 
\end{figure*}

Figure~\ref{figcheckergame} 
shows the geometry and algebra behind this discussion. 
The top shows the geometry and the bottom gives the checker description. 
It also shows the matrices parameterizing the two-dimensional families 
of lines in each case.   
The parameterization is
explicitly described in \S\ref{five}.
%%%%%%%%%%%%%%%%%%%%%%%%%%%%%%%%%%%%%%%%%%%%%%%%%%%%%%%%%%%%%%%%%%%%%%%%%%%%%%%% 
%%%%%%%%%%%%%%%%%%%%%%%%%%%%%%%%%%%%%%%%%%%%%%%%%%%%%%%%%%%%%%%%%% 
 
This single example is sufficient to understand the general case.  The 
initial position of the red checkers is as follows.  The 
intersection of the $k$-plane with the moving flag $M$ determines 
the rows of the red checkers, and the intersection with the fixed  
flag $F$ determines 
their columns, and they are arranged from southwest to northeast.  The 
movement of the moving flag in arbitrary dimension is analogous to the 
specific case described here, and is described by a sequence of moves 
of black checkers.  The movement of the black checkers determines the 
movement of the red checkers (see \S\ref{checkers}), and at 
each stage, there are one or two choices.  When there are two choices, 
the underlying geometry is essentially the same as in the 
example above.     When there is one choice, often the underlying 
geometry does not change, but the parameterization changes.

%%%%%%%%%%%%%%%%%%%%%%%%%%%%%%%%%%%%%%%%%%%%%%%%%%%%%%%%%%%%%%%%%%%%% 
% 
\section{The Littlewood-Richardson Homotopy} 
\label{six}
We first explain how the geometric  Littlewood-Richardson rule gives equations and 
homotopies for solving Schubert problems, and then illustrate that with two specific 
examples coming from the problem of four lines. 
 
In the geometric Littlewood-Richardson rule the intersection  
$\Omega_\omega(F)\cap\Omega_\tau(M)$ breaks into components which eventually become 
Schubert varieties $\Omega_\sigma(F)$ as the moving flag $M$ specializes to coincide with 
the fixed flag $F$. 
At each stage, the components correspond to checkerboards.

A checkerboard encodes the relative positions of the fixed and moving flags as well as
a representation $X$, called a localization pattern, of the general element in the
corresponding component.  
Specifically, $X$ is a $n\times k$ matrix whose entries are either $0$, $1$, or 
indeterminates such that the $n\times k$ matrix 
$MX$ is a general point in that component. 
In \S\ref{five} we explain how to obtain a localization pattern from its checkerboard. 
 
Given a Schubert problem, 
 \begin{equation}\label{Eq:SchubertProblem} 
   \Omega_\omega(F)\cap\Omega_\tau(M)\,\cap \, 
   \Omega_{\rho_1}(N^1)\cap\dotsb\cap\Omega_{\rho_s}(N^s), 
 \end{equation} 
the intersection of the last $s$ Schubert varieties is expressed 
as rank conditions on (minors of) 
matrices $[Y|N^i_j]$~\eqref{Eq:rank_condition}, where $Y$ is a general 
$n\times k$ matrix representing a general $k$-plane. 
Write this system of minors succinctly as $P(Y)=0$. 
When $X$ is a localization pattern for a checkerboard in the 
degeneration of $\Omega_\omega(F)\cap\Omega_\tau(M)$ (as in \S \ref{five}), the points in  
the corresponding component that also 
lie in the last $s$ Schubert varieties in~\eqref{Eq:SchubertProblem} 
are the solutions to the system $P(MX)=0$. 
 
Reversing the specialization of the flags $F$ and $M$ is the generalization sequence.
Between adjacent stages $i$ and $i{+}1$ of the generalization sequence, 
the moving flag is $M(t)$ for $t\in[0,1]$.
Then the homotopy connecting these stages is  
\begin{equation}
  P(M(t)X)=0 
\end{equation}
for $t\in[0,1]$. 
When $t=0$, we are in stage $i$ and when $t=1$, we are in stage $i{+}1$. 
The generalization of the moving flag is described in more detail in \S\ref{S:generalizing}.
We explain how the red checkers move in \S \ref{checkers}, and then how the localization patterns
for different stages fit together. 

We illustrate this with some examples 
from Figure~\ref{figcheckergame}.
For $X \in \Omega_{[2~4]}(F)\cap \Omega_{[2~4]}(M)$, we have
\begin{equation} 
   X = \left[ 
     \begin{array}{cc} 
        x_{11} & 0 \\ 
          1    & 0 \\ 
          0 & x_{32} \\ 
          0 &   1    \\ 
     \end{array} 
   \right] 
   \quad
   \begin{array}{l} 
      F = [ \bfe_1 , \bfe_2 , \bfe_3 , \bfe_4 ] \\
      M = [ \bfe_4, \bfe_3, \bfe_2, \bfe_1] \\ 
      {\rm for~any~} x_{11} {\rm ~and~} x_{32}: \\ 
      ~~ {\rm dim}(X \cap \langle \bfe_1 , \bfe_2 \rangle) = 1, \\ 
      ~~ {\rm dim}(X \cap \langle \bfe_4 , \bfe_3 \rangle) = 1, \\ 
      ~~ {\rm dim}(X \cap \langle \bfe_1 , \bfe_2 , \bfe_3 , \bfe_4 \rangle)
        = \makebox[0.001in][l]{$2$.}
   \end{array} 
\end{equation} 
 
%Localization patterns are encoded by red checkers, 
%see Figure~\ref{figwhitecheckers}. 
 %
%\begin{figure}[hbt] 
%\[ 
%  \includegraphics{stage0.lines.eps} 
%   \begin{picture}(64,44) 
%    \put(20,4){\begin{picture}(36,36) 
%     \put(0,27){$*$} \put(9,27){\small$0$} 
%     \put(0,18){\small$1$}   \put(9,18){\small$0$} 
%     \put(0, 9){\small$0$}   \put(9, 9){$*$} 
%     \put(0, 0){\small$0$}   \put(9, 0){\small$1$} 
%    \end{picture}} 
%    \end{picture} 
%\] 
%\caption{Red checkers encode localization patterns for the solutions.} 
%\label{figwhitecheckers} 
%\end{figure} 
 
\begin{figure}[hbt] 
\[ 
   \begin{picture}(385,80)(-95,20)
    \put(-95,10){\includegraphics[height=3.5cm]{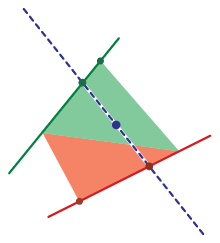}} 
    \put(195,10){\includegraphics[height=3.5cm]{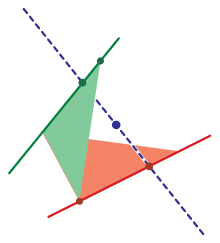}} 
 
  \put(10,41){ 
    \begin{picture}(36,36) 
     \put(0,27){$*$} \put(9,27){\small$0$} 
     \put(0,18){\small$1$}   \put(9,18){\small$0$} 
     \put(0, 9){\small$0$}   \put(9, 9){$*$} 
     \put(0, 0){\small$0$}   \put(9, 0){\small$1$} 
    \end{picture}} 
    \put(36,34){\includegraphics{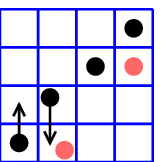}} 
%    \put(33,-9){4 3 2 1} 
    \put(95,52){\Large$\rightarrow$} 
 
  \put(175,41){ 
    \begin{picture}(36,36) 
     \put(0,27){$*$} \put(9,27){\small$0$} 
     \put(0,18){\small$1$}   \put(9,18){\small$0$} 
     \put(0, 9){\small$0$}   \put(9, 9){$*$} 
     \put(0, 0){\small$0$}   \put(9, 0){\small$1$} 
    \end{picture}} 
    \put(125,34){\includegraphics{stage1.eps}} 
%    \put(143,-9){4 3 1 2} 
 
   \end{picture}
\] 
\caption{No homotopy, only change of coordinates.}\vspace{-5pt}
\label{fignohom} 
\end{figure}

In the first stage of Figure~\ref{figcheckergame}, 
the plane in the moving flag rotates about its line  
until it meets the fixed point. 
As the line in the moving flag does not move, there is no homotopy, 
only a change of coordinates, as illustrated in Figure~\ref{fignohom} 
%%%%%%%%%%%%%%%%%%%%%%%%%%%%%%%%%%%%%%%%%%%%%%%%%%%%%%%%%%%%%%%%% 
for a line meeting two lines and a fixed point in three-space, and as
discussed at the end of \S \ref{checkers}. 
The corresponding coordinate transformation is:
\begin{eqnarray} 
   \left[ 
     \begin{array}{rrrr} 
        1 & 0 & 0 & 0 \\ 
        0 & 1 & 0 & 0 \\ 
        0 & 0 & 0 & 1 \\ 
        0 & 0 & 1 & -1 \\ 
     \end{array} 
  \right] 
  \left[ 
     \begin{array}{cc} 
        x_{11} & \!\!\!\!0 \\ 
          1 & \!\!\!\!0 \\ 
          0 & \!\!\!\!x_{32} \\ 
          0 & \!\!\!\!1 \\ 
     \end{array} 
  \right]  
  = 
  \left[ 
     \begin{array}{cc} 
        x_{11} & 0 \\ 
          1 & 0 \\ 
          0 & 1 \\ 
          0 & x_{32}{-} 1 \\ 
     \end{array} 
  \right] 
  \equiv 
  \left[ 
     \begin{array}{cc} 
        x_{11} & 0 \\ 
          1 & 0 \\ 
          0 & 1/(x_{32}{-}1) \\ 
          0 & 1 \\ 
     \end{array} 
  \right]. 
\end{eqnarray} 
 
\begin{figure}[hbt] 
\[ 
   \begin{picture}(385,80)(-95,20)
%    \put(  0,80){\includegraphics[height=3.5cm]{point1.eps}} 
    \put(-95,10){\includegraphics[height=3.5cm]{point1.eps}} 
    \put(195,10){\includegraphics[height=3.5cm]{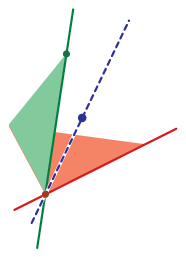}} 
 
  \put(10,41){ 
    \begin{picture}(36,36) 
     \put(0,27){\small$*$} \put(9,27){\small$0$} 
     \put(0,18){\small$1$}   \put(9,18){\small$0$} 
     \put(0, 9){\small$0$}   \put(9, 9){\small$*$} 
     \put(0, 0){\small$0$}   \put(9, 0){\small$1$} 
    \end{picture}} 
    \put(36,34){\includegraphics{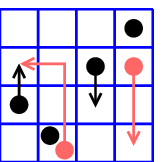}} 
%    \put(33,-9){4 3 1 2} 
    \put(95,52){\Large$\rightarrow$} 
 
  \put(175,41){ 
    \begin{picture}(36,36) 
     \put(0,27){\small$0$} \put(9,27){\small$*$} 
     \put(0,18){\small$1$}   \put(9,18){\small$0$} 
     \put(0, 9){\small$0$}   \put(9, 9){\small$*$} 
     \put(0, 0){\small$0$}   \put(9, 0){\small$1$} 
    \end{picture}} 
    \put(125,34){\includegraphics{stage2.1.eps}} 
%    \put(143,-9){4 1 3 2} 
 
   \end{picture}
\] 
\caption{Homotopy, as red checkers swap rows.} 
\label{figlrhomswap} 
\end{figure}

When the red checkers swap rows, we use a homotopy, 
shown in Figure~\ref{figlrhomswap}, also for the case 
of a line meeting two lines and a fixed point. 
This homotopy has coordinates 
\begin{equation} 
X(t) = 
 \left[ 
   \begin{array}{lc} 
     x_{12} t &  x_{12} \\ 
     x_{32}   &  0      \\ 
     x_{32} t &  x_{32} \\ 
     \multicolumn{1}{c}{0}      &  1 
   \end{array} 
 \right] \ .
\end{equation} 
At $t=0$ we see that $X(0)$ fits the pattern on the right in Figure~\ref{figlrhomswap}, 
while at $t = 1$ a coordinate change 
brings $X(1)$ into the pattern on the left. 
With linear combinations of the two columns we find generators for
the line that fit the columns of the pattern.
 
%%%%%%%%%%%%%%%%%%%%%%%%%%%%%%%%%%%%%%%%%%%%%%%%%%%%%%%%%%%%%%%%%%%%% 
% 
\section{Localization Patterns} 
 \label{five}
We describe coordinates for each component corresponding to 
a checkerboard:  given two flags $M$ and $F$ 
in relative position described by the black checkers, it is the space 
of $k$-planes meeting $M$ and $F$ in the manner specified by the 
positions of the red checkers. 
The black checkers correspond to a basis of both $F$ and $M$.
Each red checker is a basis element for the $k$-plane and it lies in the space spanned
by the black checkers weakly to its northwest.

\begin{figure}[hbt]
\[
  \raisebox{-77pt}{\begin{picture}(154,154)
   \put(0,0){\includegraphics{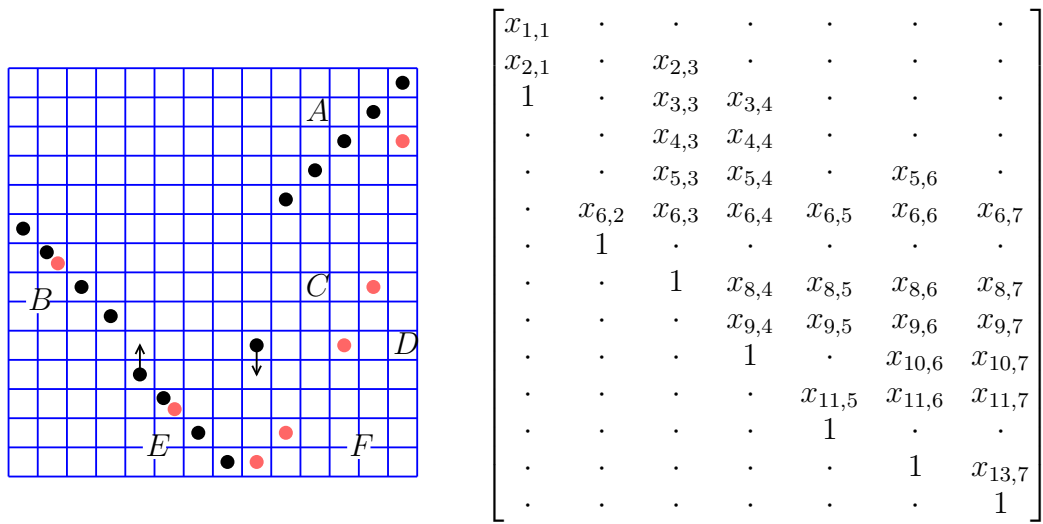}}\
   \put(112.5,134.5){$A$}
   \put(  7.5, 63){$B$}
   \put(112.5, 68.5){$C$}
   \put(145.5, 46.5){$D$}
   \put( 52,  8){$E$}
   \put(129,  8){$F$}
  \end{picture}}\qquad
\left[
\begin{matrix}
 x_{1,1}&\cdot &\cdot &\cdot&\cdot &\cdot  &\cdot\\ 
 x_{2,1}&\cdot &x_{2,3} &\cdot&\cdot &\cdot &\cdot\\ 
  1    &\cdot &x_{3,3}&x_{3,4}&\cdot &\cdot &\cdot\\ 
 \cdot &\cdot &x_{4,3}&x_{4,4}&\cdot &\cdot  &\cdot\\ 
 \cdot &\cdot &x_{5,3}&x_{5,4}&\cdot &  x_{5,6}&\cdot\\ 
 \cdot &x_{6,2}&x_{6,3} &x_{6,4}&x_{6,5}&x_{6,6}  &x_{6,7}\\ 
 \cdot &  1   &\cdot &\cdot &\cdot &\cdot  &\cdot     \\ 
 \cdot &\cdot &1     &x_{8,4}&x_{8,5}&x_{8,6}  &x_{8,7}\\ 
 \cdot &\cdot &\cdot &x_{9,4}&x_{9,5}&x_{9,6}  &x_{9,7}\\ 
 \cdot &\cdot &\cdot & 1    &\cdot &x_{10,6} &x_{10,7}\\ 
 \cdot &\cdot &\cdot &\cdot &x_{11,5}&x_{11,6}&x_{11,7}\\ 
 \cdot &\cdot &\cdot &\cdot & 1     &\cdot &\cdot\\ 
 \cdot &\cdot &\cdot &\cdot &\cdot  &1     &x_{13,7}\\ 
 \cdot &\cdot &\cdot &\cdot &\cdot  &\cdot  &1  
\end{matrix}\right]
\]
\caption{Coordinates corresponding to a checkerboard.
         Entries $\cdot$ in the coordinate matrix are $0$.} 
\label{t:ttt}
%\end{center}
\end{figure}

While special cases were shown in Figure~\ref{figcheckergame}, we
illustrate the general case with an example.  
In the checkerboard of Figure~\ref{t:ttt}, 
one black checker (in row $D$) is descending.   
Red checkers are distributed along the sorted black checkers 
(regions $B$ and $E$), as well as in the pre-sorted region 
(regions $A$, $C$, $D$, and $F$); in the latter region,
they are distributed from the southwest to the northeast as shown.  
The corresponding localization pattern, which is expressed 
with respect to the basis of $M$, is shown in Figure~\ref{t:ttt}.  

We discuss the linking of localization patterns between stages 
after we describe the movement of checkers in \S \ref{checkers}. 

%%%%%%%%%%%%%%%%%%%%%%%%%%%%%%%%%%%%%%%%%%%%%%%%%%%%%%%%%%%%%%%%%% 
\section{Generalizing the moving Flag} 
\label{S:generalizing}
 
Underlying the geometric Littlewood-Richardson rule is the sequence of specializations 
(analogous to Figure~\ref{figmovingflag}) in which the moving flag $M$ successively 
moves to coincide with the fixed flag $F$. 
Reversing this gives the generalization sequence in which $M$ emerges from $F$.  
 
The generalization of the moving flag $M$ is  as follows. 
Throughout, the fixed flag is 
\begin{equation}
F = \left\{ \langle \bfe_1 \rangle \subset \langle \bfe_1, \bfe_2 \rangle 
  \subset  \langle \bfe_1, \bfe_2, \bfe_3 \rangle \cdots \right\}.
\end{equation}
Initially, the moving flag $M$ coincides with $F$.   
We let $\bfm'_i(t)$ describe the 
vectors during the generalization ($t=0$ corresponds to the 
specialized case, and $t=1$ corresponds to the generalized  case), and  
$\bfm''_i$ describe the vectors after the generalization.   
At time $t$,  
\begin{equation}
  M(t) = \left\{ \langle \bfm'_1(t)  \rangle \subset 
   \langle \bfm'_1(t), \bfm'_2(t) \rangle 
  \subset \cdots \right\}.
\end{equation}
In the checker diagram, at each stage the black checkers in rows $r$ and $r+1$ swap rows,
for some  $r$.  
Set 
\begin{equation}
    { \bfm_i = \bfm'_i(t) = \bfm''_i \quad \quad 
    \mbox{for ${ i \neq r, r+1}$}.} 
\end{equation}
These different notations for the same vector keep track of  
whether we are talking about $t=0$, general $t$, or $t=1$. 
\begin{align} 
 \bfm_r       &=  \bfm''_{r+1} (=\bfm'_r(0) = \bfm'_{r+1}(0)),   \label{star1}\\ 
 \bfm_{r+1}    &= \bfm''_{r+1} -  \bfm''_r  (=  \bfm'_{r+1}(0)- \bfm'_r(0) ). \\ 
 \bfm'_r(t)   &= t \bfm''_r + (1-t) \bfm''_{r+1} = \bfm''_{r+1} - t \bfm_{r+1},  \\  
 \bfm'_{i}(t) &= \bfm''_{i}  \quad \mbox{ for all other $i$}. \label{star2}
\end{align} 
Thus $\bfm'_i(1) = \bfm''_i$ for all $i$.   
 
It is convenient to describe the homotopy in terms of matrices.  
Here are the generalizing moves from Figure~\ref{figmovingflag}.   
{\small
\begin{eqnarray} 
\quad F = \left[ 
    \begin{array}{cccc} 
       1 & 0 & 0 & 0 \\  
       0 & 1 & 0 & 0 \\  
       0 & 0 & 1 & 0\\  
       0 & 0 & 0 & 1  
    \end{array} 
  \right] 
  \rightarrow 
  \left[ 
    \begin{array}{cccc} 
       1 & 0 & 0 & 0 \\  
       0 & 1 & 0 & 0 \\  
       0 & 0 & \gamma_{31} & 1 \\  
       0 & 0 & 1 & 0   
    \end{array} 
  \right]
%\end{eqnarray} 
%\begin{eqnarray} 
\rightarrow 
  \left[ 
    \begin{array}{cccc} 
       1 & 0 & 0 & 0 \\  
       0 & \gamma_{21} & 1 & 0 \\  
       0 & \gamma_{31} & 0 & 1 \\  
       0 & 1 & 0 & 0 \\  
    \end{array} 
  \right]  
\rightarrow 
  \left[ 
    \begin{array}{cccc} 
       \gamma_{11} & 1 & 0 & 0 \\  
       \gamma_{21} & 0 & 1 & 0 \\  
       \gamma_{31} & 0 & 0 & 1 \\  
       1 & 0 & 0 & 0 \\  
    \end{array} 
  \right] \\ 
%\end{eqnarray} 
%\begin{eqnarray} 
\rightarrow 
  \left[ 
    \begin{array}{cccc} 
       \gamma_{11} & 1 & 0 & 0 \\  
       \gamma_{21} & 0 & \gamma_{22} & 1 \\  
       \gamma_{31} & 0 & 1 & 0 \\  
       1 & 0 & 0 & 0 \\  
    \end{array} 
  \right]  
\rightarrow \label{Eq:nineteen}
  \left[ 
    \begin{array}{cccc} 
       \gamma_{11} & \gamma_{12} & 1 & 0 \\  
       \gamma_{21} & \gamma_{22} & 0 & 1 \\  
       \gamma_{31} & 1 & 0 & 0 \\  
       1 & 0 & 0 & 0 \\  
    \end{array} 
  \right] 
%\end{eqnarray} 
%\begin{eqnarray} 
 \rightarrow  
  \left[ 
    \begin{array}{cccc} 
       \gamma_{11} & \gamma_{12} & \gamma_{13} & 1 \\  
       \gamma_{21} & \gamma_{22} & 1 & 0 \\  
       \gamma_{31} & 1 & 0 & 0 \\  
       1 & 0 & 0 & 0 \\  
    \end{array} 
  \right]. 
\end{eqnarray} }
Here, $\gamma_{ij}$ are general complex numbers. 
For example, the second matrix in~\eqref{Eq:nineteen} corresponds to stage $1$, and 
we see that the moving plane, (the projectivization of)  the span of
the first three columns,
indeed contains the fixed point, as $e_1$ is in the span of those
three column vectors, in agreement with Figure~\ref{figmovingflag}.
 
The arrows represent the movement of the flag $M$, 
which we parametrize using our homotopy parameter $t\in[0,1]$. 
For example, the next to last deformation is 
\begin{eqnarray} 
  \left[ 
    \begin{array}{cccc} 
       \gamma_{11} & 1 & 0 & 0 \\  
       \gamma_{21} & 0 & \gamma_{22} & 1 \\  
       \gamma_{31} & 0 & 1 & 0 \\  
       1 & 0 & 0 & 0 \\  
    \end{array} 
  \right] 
  \left[ 
    \begin{array}{cccc} 
       1 & 0 & 0 & 0 \\  
       0 & \gamma_{12} t & 1 & 0 \\  
       0 & 1 & 0 & 0\\  
       0 & 0 & 0 & 1 \\  
    \end{array} 
  \right] \\ 
  =  
  \left[ 
    \begin{array}{cccc} 
       \gamma_{11} & \gamma_{12} t & 1 & 0 \\  
       \gamma_{21} & \gamma_{22} & 0 & 1 \\  
       \gamma_{31} & 1 & 0 & 0 \\  
       1 & 0 & 0 & 0 \\  
    \end{array} 
  \right]\ =:M(t).\label{Eq:M(t)} 
\end{eqnarray} 
The gradual introduction of the random constants $\gamma_{ij}$ in the moving flag 
is the analog here of the gamma trick~\cite{SW05} to ensure the 
regularity of the solution paths.  By this gamma trick,  
for all~$t$, except for a finite number of choices of $\gamma_{ij}$,  
the solution paths contain only regular points. 
 
The Littlewood-Richardson homotopies operate on randomly generated 
complex flags.  To move to flags with specific coordinates,  
we use coefficient-parameter~\cite{MS89} or cheater homotopies~\cite{LSY89}. 
  
%%%%%%%%%%%%%%%%%%%%%%%%%%%%%%%%%%%%%%%%%%%%%%%%%%%%%%%%%%%%%%%%%%%%% 
% 
\section{Movement of red checkers} 
\label{checkers}

\begin{figure}[hbt]
\[
  \begin{picture}(224,74)(-20,-2)
   \put(80,-2){\includegraphics{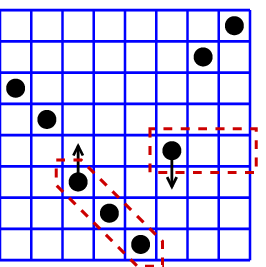}}
   \put(158,29){critical row}
   \put( -20,11){critical diagonal}
   \put( 68,13.5){\vector(1,0){28}}
  \end{picture}
\]
\caption{Critical row and critical diagonal.} 
\label{F:risefall}
\end{figure}

In the geometric Littlewood-Richardson rule, the black checkers start out 
on the anti-diagonal,
and a bubble sort is performed which moves them to the diagonal.
This is indicated in Figures~\ref{figmovingflag},~\ref{figcheckergame}, 
and~\ref{t:ttt}.
In each of the $\binom{n}{2}$ steps, one black checker descends and another 
rises as in Figure~\ref{F:risefall}.
%%%%%%%%%%%%%%%%%%%%%%%%%%%%%%%%%%%%%%%%%%%%%%%%%%%%%%%%%%%%%%%%%%%%%%%%%%%%%
%%%%%%%%%%%%%%%%%%%%%%%%%%%%%%%%%%%%%%%%%%%%%%%%%%%%%%%%%%%%%%%%%%%%%%%%%%%%%
The descending checker is in the critical row and the ascending checker is
at the top left of the critical diagonal.

To resolve the intersection $\Omega_\omega(F)\cap\Omega_\tau(M)$, 
we initially place red checkers as follows.  The 
intersection of the $k$-plane with the moving flag $M$ determines 
the rows of the red checkers, and the intersection with the fixed  
flag $F$ determines 
their columns, and they are arranged from southwest to northeast. 
As the black checkers move, they induce a motion of the red checkers.
There will be nine cases to consider.
In eight, the motion is determined, while in the ninth case there are 
sometimes two choices as in Figure~\ref{figcheckergame}.

The cases are determined by the answers to two questions, 
each of which has three
answers. 
\begin{enumerate}
 \item Where is the top red checker in the critical diagonal?
  \begin{itemize}
   \item[$(a)$] In the rising checker's square.
   \item[$(b)$] Elsewhere in the critical diagonal.
   \item[$(c)$] There is no red checker in the critical diagonal.
  \end{itemize}
 \item Where is the red checker in the critical row?
  \begin{itemize}
   \item[$(\alpha)$] In the descending checker's square.
   \item[$(\beta)$] Elsewhere in the critical row.
   \item[$(\gamma)$] There is no red checker in the critical row.
  \end{itemize}
\end{enumerate}

Table~\ref{t:move} shows the movement of the checkers in these nine cases.
The rows correspond to the answers to the first question and the columns to the answers of 
the second question.
Only the relevant part of each checkerboard is shown.

%%%%%%%%%%%%%%%%%%%%%%%%%%%%%%%%%%%%%%%%%%%%%%%%%%%%%%%%%%%%%%%%%%%%%%%%%%%%%
\begin{table}[hbt]
\begin{center}
 \begin{tabular}{|c|c|c|c|}\hline
 & $\alpha$&$\beta$&$\gamma$\\\hline
 \raisebox{15pt}{$a$}
 &\raisebox{5.5pt}{\includegraphics{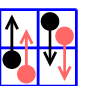}}& \includegraphics{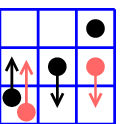}\rule{0pt}{39pt}&
 \raisebox{5.5pt}{\includegraphics{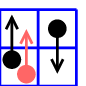}}\\\hline
 \raisebox{19.5pt}{$b$}
 & \raisebox{5.5pt}{\includegraphics{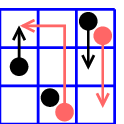}} & \includegraphics{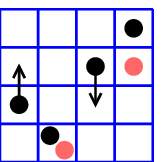}\rule{0pt}{50pt}
   \raisebox{18pt}{\ or\ }  \includegraphics{22b.eps} &
   \raisebox{5.5pt}{\includegraphics{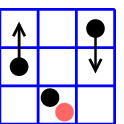}} \\\hline
 \raisebox{15pt}{$c$}&
 \raisebox{5.5pt}{\includegraphics{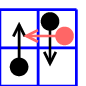}} & \includegraphics{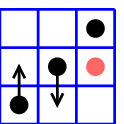}&
 \raisebox{5.5pt}{\includegraphics{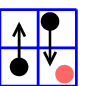}}\rule{0pt}{39pt}\\\hline
\end{tabular}
\end{center}
\caption{Movement of red checkers.} 
\label{t:move}
\end{table}
%%%%%%%%%%%%%%%%%%%%%%%%%%%%%%%%%%%%%%%%%%%%%%%%%%%%%%%%%%%%%%%%%%%%%%%%%%%%%

\begin{figure}[hbt]
\[
 \begin{picture}(270,47.5)(-40,5)
   \put(70,0){\includegraphics{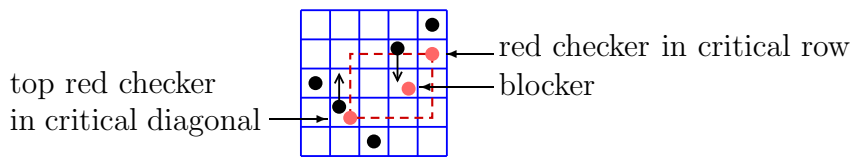}}
   \put(145,38){red checker in critical row}
   \put(143,38.5){\vector(-1,0){18}}

   \put(-40,17){\parbox{95pt}{top red checker\newline  in critical diagonal}}
   \put(58,14){\vector(1,0){22}}

   \put(145,23){blocker}\put(143,26){\vector(-1,0){28}}
 \end{picture}
\]
\caption{a blocker.} 
\label{F:blocker}
\end{figure}

In case $(b,\beta)$ there are two possibilities, 
which can both occur---this is when a
component breaks into two components in the geometric 
Littlewood-Richardson rule.
The second of these (where the red checkers swap rows) only occurs 
if there are no other
red checkers in the rectangle between the two, which we call blockers.
Figure~\ref{F:blocker} shows a blocker.
%%%%%%%%%%%%%%%%%%%%%%%%%%%%%%%%%%%%%%%%%%%%%%%%%%%%%%%%%%%%%%%%%%%%%%%%%%%%%
%%%%%%%%%%%%%%%%%%%%%%%%%%%%%%%%%%%%%%%%%%%%%%%%%%%%%%%%%%%%%%%%%%%%%%%%%%%%%

To track solutions to Schubert problems between adjacent stages
in the generalization sequence, 
we need uniform coordinates corresponding to two adjacent 
diagrams---for example, two boards connected by an arrow in 
Figure~\ref{figcheckergame}.  We can then track solutions from one 
board to the more generalized board.  There are three cases 
to consider.   
 
In trivial cases, such as the first arrow in 
Figure~\ref{figcheckergame}, which is case $(c,\gamma)$ of Table~\ref{t:move}, the
coordinates do not change because the underlying geometry is constant.  
 
We describe one of the nontrivial examples of the coordinates linking 
two stages, that of the lower arrow between stage 1 and stage 2 in 
Figure~\ref{figcheckergame} (left case of $(b,\beta)$). 
We follow the vector corresponding to the red checker in the bottom 
row.  Throughout the degeneration (as $t$ goes from $0$ to $1$), we 
write its vector 
as $\bfm_4 + x \bfm'_2(t)$.  As $\bfm'_2(1) = \bfm''_2$ and $\bfm'_2(0) = \bfm''_3$ 
(see \eqref{star1}--\eqref{star2} of \S \ref{S:generalizing} with $r=2$), 
we see that the reason for the change of row of 
the $*$ in the matrix in Figure~\ref{figcheckergame} is just a 
renaming of the variable.   
 
The third case, where the two red checkers swap rows, is 
more subtle, and an example was described at the end of \S \ref{six}. 

%%%%%%%%%%%%%%%%%%%%%%%%%%%%%%%%%%%%%%%%%%%%%%%%%%%%%%%%%%%%%%%%%%% 
\section{Solving Schubert Problems}
\label{S:solving} 
 
The global structure of the Littlewood-Richardson homotopy is encoded by a graded poset.
This records the branching of Schubert varieties that occur in when running the geometric
Littlewood-Richardson rule through successive specializations of their defining flags,
equivalently, moving checkers as in \S \ref{checkers}.

We construct the poset for a a Schubert problem
\begin{equation}
   \Omega_{\omega^1}(F^1)\cap
   \Omega_{\omega^2}(F^2)\cap\dotsb\cap
   \Omega_{\omega^s}(F^s).
\end{equation}
First, use the geometric Littlewood-Richardson rule to resolve 
the first intersection
 \begin{equation}\label{Eq:resolve}
    \Omega_{\omega^1}(F^1)\cap
   \Omega_{\omega^2}(F^2) \sim
   \sum_{\sigma}c^{\sigma}_{\omega^1,\omega^2} \Omega_\sigma(F^1).
 \end{equation}
The top of the poset is the bracket $\omega^1$,
which branches to those brackets $\sigma$ appearing in the sum.
The edge $\omega^1\to\sigma$ occurs with 
multiplicity $c^{\sigma}_{\omega^1,\omega^2}$.
Geometrically, we have the disjunction of Schubert probems
\begin{equation}
   \Bigl(\sum_{\sigma}c^{\sigma}_{\omega^1,\omega^2} \Omega_\sigma(F^1)\Bigr)
   \cap \Omega_{\omega^3}(F^3)\cap\dotsb\cap
   \Omega_{\omega^s}(F^s),
\end{equation}
and we resolve each $\Omega_\sigma(F^1)\cap\Omega_{\omega^3}(F^3)$ 
with the geometric Littlewood-Richardson rule, further building the poset,
and continue in this fashion.

The penultimate stage has the form
\begin{equation}
  \Bigl(\sum_\sigma C^\sigma \Omega_\sigma(F^1)\Bigr)
   \cap \Omega_{\omega^s}(F^s),
\end{equation}
where $C^\sigma$ are the multiplicities.
This is resolved via~\eqref{Eq:double_intersection}, so the only term in the 
sum which contributes is when $\sigma^\vee=\omega^s$, and the 
final Schubert variety is $\Omega_{[1~2~\dotsb~k]}(F^1)$.

The global structure of the Littlewood-Richardson homotopy is 
to begin with the solution
$\Omega_{[1~2~\dotsb~k]}(F^1)$ at the bottom of our poset, 
and continue this solution
along homotopies corresponding to the edges of the poset.
Each edge is a sequence of $\binom{n}{2}$ homotopies or coordinate 
changes corresponding to running the geometric Littlewood-Richardson rule backwards, as
explained in \S \ref{five}, \S\ref{S:generalizing}, and \S\ref{checkers}. 
In this way, we iteratively build solutions to the Schubert-type
problems corresponding to the nodes of this poset.

For example, suppose that we have the Schubert problem $[2~4~6] [2~5~6]^3$.
This is resolved in the geometric Littlewood-Richardson rule as 
\begin{align} 
[2\:4\:6] [2\:5\:6]^3
 &=  (1[2\:3\:5]+1[1\:4\:5]+1[1\:3\:6]) [2\:5\:6]^2 \\ 
 &= (2[1\:3\:4]+2[1\:2\:5]) [2\:5\:6] \\ 
 &= 2[1\:2\:3]. 
\end{align} 

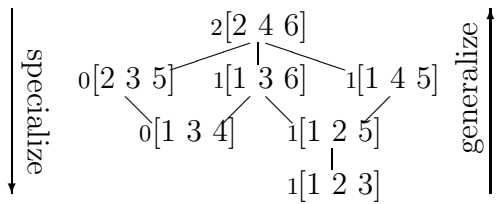
\begin{figure}[hbt] 
\begin{center} 
\begin{picture}(180,70)
\put(0,70){\vector(0,-1){70}} 
\put(6,55){\begin{rotate}{-90}specialize\end{rotate}}
%\put(1,64){s}  \put(190,64){g} 
%\put(0,57){p}  \put(190,57){e} 
%\put(0,50){e}  \put(190,50){n} 
%\put(0,43){c}  \put(190,43){e} 
%\put(1,36){i}  \put(190,36){r} 
%\put(0,30){a}  \put(190,30){a} 
%\put(1,23){l}  \put(191,23){l} 
%\put(1,15){i}  \put(191,15){i} 
%\put(0,9){z}   \put(190,9){z} 
%\put(0,2){e}   \put(190,2){e} 

\put(176,15){\begin{rotate}{90}generalize\end{rotate}}
\put(181,0){\vector(0,+1){70}} 
\put(75,60){\scriptsize 2} \put(80,60){[2~4~6]} 
\put(90,57){\line(-3,-1){30}} 
\put(93,57){\line(0,-1){8}} 
\put(96,57){\line(+3,-1){30}} 
\put(25,40){\scriptsize 0}   \put(30,40){[2~3~5]}  
\put(76,40){\scriptsize 1}  \put(80,40){[1~3~6]}  
\put(126,40){\scriptsize 1} \put(130,40){[1~4~5]} 
\put(43,37){\line(1,-1){10}} 
\put(90,37){\line(-1,-1){10}} 
\put(96,37){\line(1,-1){10}} 
\put(143,37){\line(-1,-1){10}} 
\put(48,20){\scriptsize 0}  \put(53,20){[1~3~4]}  
\put(104,20){\scriptsize 1} \put(108,20){[1~2~5]} 
\put(121,17){\line(0,-1){8}} 
\put(104,0){\scriptsize 1}  \put(108,0){[1~2~3]} 
\end{picture} 
\caption{Poset to resolve $[2~4~6] [2~5~6]^3$.} 
\label{figposet} 
\end{center} 
\end{figure}

The poset corresponding to  
the Littlewood-Richardson homotopies is shown in Figure~\ref{figposet}. 
%%%%%%%%%%%%%%%%%%%%%%%%%%%%%%%%%%%%%%%%%%%%%%%%%%%%%%%%%%%%
%%%%%%%%%%%%%%%%%%%%%%%%%%%%%%%%%%%%%%%%%%%%%%%%%%
The multiplicities in front of the brackets are the number of solutions
tracked to the given Schubert variety.
 
The number of solution paths is one of the three factors that 
determine the cost of the homotopies.  Another factor is the complexity 
of the polynomials that express the intersection conditions. 
The current implementation performs a Laplace expansion on the minors 
to elaborate all conditions~\eqref{Eq:Schubert_condition}. 
Locally, for use during path following, an overdetermined system of $p$ 
equations in $q$ unknowns is multiplied with a $q$-by-$p$ matrix of 
randomly generated complex coefficients to obtain square linear systems 
in the application of Newton's method.  The third factor in the cost 
lies in the complex and real geometry of the solution paths.   
In practice it turns out that 
solving a generic complex instance with the Pieri homotopies  
is in general always faster than running a cheater homotopy  
using the solutions of a generic complex instance as start solutions 
to solve a generic real instance. 
This experience also applies to solving general Schubert problems 
with Littlewood-Richardson homotopies. 
 
%%%%%%%%%%%%%%%%%%%%%%%%%%%%%%%%%%%%%%%%%%%%%%%%%%%%%%%%%%
%
\section{Computational Experiments}
\label{S:computation} 
 
Littlewood-Richardson homotopies are available in PHCpack~\cite{Ver99} 
since release 2.3.46.  Release 2.3.52 contains {\tt LRhomotopies.m2}, 
an interface to solve Schubert problems in Macaulay~2~\cite{macaulay2}.  
Via {\tt phc -e} option \#4 we resolve intersection conditions 
and Littlewood-Richardson homotopies are available via option~\#5.  
 
Below we list sample timings for  
solving some small Schubert problems on one core of a Mac OS X 2.2 Ghz: 
\newline $\bullet$ $[2~4]^4 = 2$ takes 5 milliseconds, 
\newline $\bullet$ $[2~4~6]^3 = 2$ takes 169 milliseconds, 
\newline $\bullet$ $[2~5~8]^2 [4~6~8] = 2$ takes 2.556 seconds, 
\newline $\bullet$ $[2~4~6~8]^2 [2~5~7~8] = 3$ takes 8.595 seconds. 
 
%\newpage 
%%%%%%%%%%%%%%%%%%%%%%%%%%%%%%%%%%%%%%%%%%%%%%%%%%%%%%%%%%

\end{document}